\documentclass[10pt,twoside]{article}
\usepackage{graphicx}
\usepackage{amsmath}
\usepackage{Latex-document}

\title{\bf  Biological Sequence Analysis\vskip 6mm}

\author{T. P. Speed\vspace*{-0.5cm}\thanks{Department of Statistics,
University of California, Berkeley, CA 94720, USA; Division of Genetics and
Bioinformatics, Walter and Eliza Hall Institute of Medical Research, VIC
3050, Australia.
E-mail: terry@stat.berkeley.edu}}
\date{\vspace{-8mm}}

\begin{document}

\maketitle

\thispagestyle{first} \setcounter{page}{97}

\begin{abstract}\vskip 3mm
This talk will review a little over a decade's research on applying
certain stochastic models to biological sequence analysis. The models
themselves have a longer history, going back over 30
years, although many novel variants have arisen since that time. The
function of the models in biological sequence analysis is to summarize
the information concerning what is known as a motif or a domain in
bioinformatics, and to provide a tool for discovering instances of that
motif or domain in a separate sequence segment.
We will introduce the motif models in stages, beginning from very simple,
non-stochastic versions, progressively becoming more complex, until we
reach modern profile HMMs for motifs. A second example will come from gene
finding using sequence data from one or two species, where generalized
HMMs or generalized pair HMMs have proved to be very effective.
\vskip 4.5mm


\end{abstract}

\vskip 12mm

\section{Introduction} \label{section 1}\setzero
\vskip-5mm \hspace{5mm}

DNA (deoxyribonucleic acid),  RNA (ribonucleic acid), and proteins are
macromolecules which are unbranched polymers
built up from smaller units. In the case of DNA these units are the 4
nucleotide residues A (adenine), C (cytosine), G (guanine) and T (thymine)
while for RNA the units are the 4 nucleotide residues A, C, G and U (uracil).
For proteins the units are the 20 amino acid residues A (alanine),
C (cysteine) D (aspartic acid), E (glutamic acid), F (phenylalanine),
G (glycine), H (histidine), I (isoleucine), K (lysine), L (leucine),
M (methionine),  N (asparagine), P (proline), Q (glutamine), R
(arginine), S (serine), T (threonine), V (valine), W (tryptophan) and Y
(tyrosine). To a considerable extent, the chemical properties of
DNA, RNA and protein molecules are encoded in the linear
sequence of these basic units: their primary structure. %%%\\

The use of statistics to study linear sequences of biomolecular units
can be descriptive or it can be predictive. A very wide range of
statistical techniques has been used in this context, and while
statistical models can be extremely useful, the underlying stochastic
mechanisms should never be taken literally. A model or method
can break down at any time without notice. Further, biological
confirmation of predictions is almost always necessary.%%%\\

The statistics of biological sequences can be global or it can be
local. For example, we might consider the global base composition of genomes:
{\em E. coli} has 25\% A, 25\% C, 25\% G, 25\% T, while
{\em P. falciparum} has 82\%A+T.  %%%\\
At the very local, the triple ATG is the near universal motif
indicating the start of translation in DNA coding sequence.
A major role of statistics in this context is to characterize
individual sequences or classes of biological sequences using
probability models, and to make use of these models to identify
them against a background of other sequences. Needless to say, the
models and the tools vary greatly in complexity. %%%\\

Extensive use is made in biological sequence analysis of the notions
of motif or domain in proteins, and site in DNA.  We shall use these
terms interchangeably to describe the recurring elements of interest to
us. It is important to note that while we focus on the sequence
characteristics of motifs, domains or sites, in practice they also
embody (biochemical) structural significance.

\section{Deterministic models} \label{section 2}
\setzero\vskip-5mm \hspace{5mm}

\begin{figure}[h]
  \centerline{\includegraphics[width=10cm]{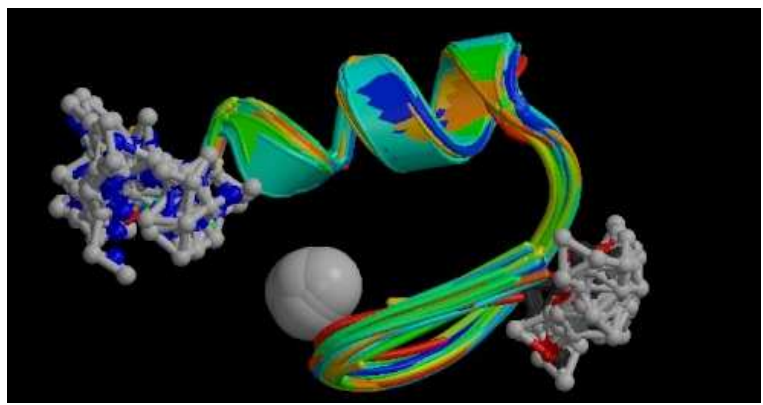}}
  \caption{A C2H2 zinc finger DNA binding domain}
  \label{fig1}
\end{figure}

The C2H2 (cysteine-cysteine histidine-histidine) zinc-finger
DNA binding domain is composed of
25-30 amino acid residues including two conserved cysteines and two
conserved histidines spaced in a particular way, with some restrictions
on the residues in between and nearby. Of course the arrangement
reflects the three-dimensional molecular structure into which the
amino-acid sequence folds, for it is the structure which has the real
biochemical significance, see Figure~\ref{fig1}, which was obtained from
{\tt http://www.rcsb.org/pdb/}.
An example of this motif is the 27-letter sequence known as 1ZNF,
this being a Protein Data Bank identifier for the structure XFIN-31 of
{\em X. laevis}. Its amino acid sequence is
\begin{verbatim}
1ZNF:  XYKCGLCERSFVEKSALSRHQRVHKNX
\end{verbatim}
Note the presence of the two $C$s separated by 2 other residues,
and the two $H$s separated by 3 other residues. Here and
elsewhere, X denotes an arbitrary amino acid residue.
A popular and useful summary description of C2H2 zinc fingers which
clearly includes our example, is the regular expression
%%%$$
\begin{equation*}
C-X(2,4)-C-X(3)-[LIVMFYWC]-X(8)-H-X(3,5)-H
\end{equation*}
%%%$$
where $X(m)$ denotes a sequence of $n$ unspecified amino acids,
while $X(m,n)$ denotes from $m$ to $n$ such, and the brackets
enclose mutually exclusive alternatives. There is a richer set of
notation for {\em regular expressions} of this kind, but for our
purposes it is enough to note that this representation is
essentially deterministic, with uncertainty included only through
mutually exclusive possibilities (e.g. length or residue) which
are not otherwise distinguished.

Simple and efficient algorithms exist for searching query
sequences of residues to find every instance of the regular
expression above. In so doing with sequence in which all
instances of the motif are known, we may identify some
sub-sequences of the query sequence which are not C2H2 zinc
finger DNA binding domains, i.e. which are false positives,  and
we may miss some sub-sequences which are C2H2 zinc fingers, i.e.
which are false negatives. Thus we have essentially deterministic
descriptions and search algorithms for the C2H2 motifs using
regular expressions. Their performance can be described by the
frequency of false positives and false negatives, equivalently,
their complements, the specificity and sensitivity of the regular
expression. We do not have space for an extensive bibliography,
so for more on regular expressions and on most of the other
concepts we introduce below, see \cite{durbin:1998}.

\section{Regular expressions can be limiting} \label{section 3}
\setzero\vskip-5mm \hspace{5mm}

Most protein binding sites are characterized by some degree
of sequence specificity, but seeking a consensus DNA sequence
is often an inadequate way to recognize their motifs.
Simply listing the alternatives seen at a position
may not be very informative, but keeping track of the
frequencies with which the different alternatives appear
can be very valuable. Thus position-specific
nucleotide or amino acid distributions came to represent
the variability in DNA or protein motif composition. This is just the
set of marginal distribution of letters at each position. Rather
than present an extensive tabulation of frequencies for our
C2H2 zinc finger example, we present a pictorial representation:
a sequence logo coming from {\tt http://blocks.fhcrc.org}. \\
%%
%%%%%
\begin{figure}[h]
\centerline{\includegraphics[width=10cm]{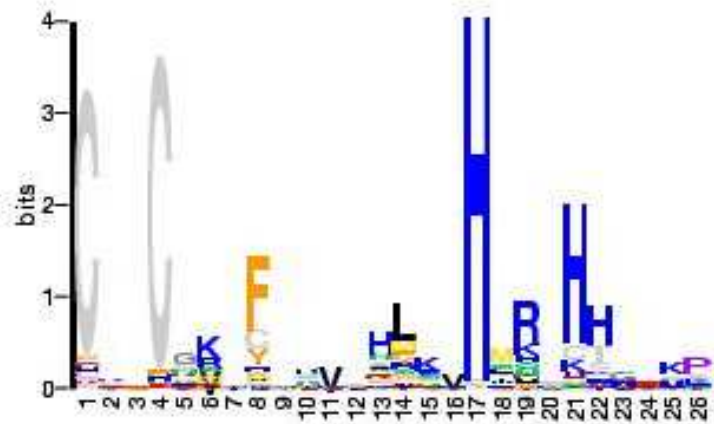}}
  \caption{Sequence logo for C2H2 zinc finger}
  \label{fig2}
\end{figure}
%%%%%
%%%Figure 2 here\\
%%

Sequence logos are
scaled representation of position-specific nucleotide or amino acid
distributions. The overall height at a given position is proportional to
information content, which is a constant minus the
entropy of the distribution at that position. The proportions of
each nucleotide or amino acid at a position are in relation to
their observed frequency at that position, with the most frequent
on top, the next most frequent below, etc. %%%\\

\section{Profiles} \label{section 4}
\setzero\vskip-5mm \hspace{5mm}

It is convenient for our present purposes to define a profile as a set
of position-specific distributions describing a motif. (Traditionally
the term has been used for the derived scores.) How would we use a set
of such distributions to search a query sequence for instances of
the motif? The answer from bioinformatics is that we {\em score} the query
sequence, and for suitably large scores, declare that a candidate
subsequence is an instance of our motif. %%%\\

There are a number of approaches for deriving profile scores,
but the easiest to explain here is this:
scores are {\em log-likelihood ratio test statistics}, for discriminating
between a probability model $M$ for the motif and a model $B$ for the
background. The model $M$ will be the direct product of
the position-specific distributions, (i.e. the independent but not
identical distribution model), while the background model $B$
will be the direct product of a set of relevant background frequencies
(i.e. the independent and identical distribution model).
Thus, if $f_{al}$ is the frequency of residue $a$ at
position $l$ of the motif, and $f_a$ background frequency of the same
residue, then the profile score assigned to residue $a$ at position
$l$ in a possible instance of the motif will be
$s_{al} = \log f_{al}/f_a$.
These scores are then summed across the positions
in the motif, and compared to a suitably defined threshold.
Note that proper setting of the threshold requires a set of data
in which all instances of the motif are known. The false positive and
false negative rate could then be determined for various thresholds,
and a suitable choice made.

We briefly discuss variants of the log-likelihood ratio scores.
In many contexts, it will matter little whether a position is
occupied by a leucine ($L$) rather than an isoleucine ($I$), as each
can evolve in time to or from the other rather more readily than from other
residues. Thus it might make sense to modify the scores to take this and
similar evolutionary patterns into account. Indeed the first use of profiles
involved scores of this kind, using the position specific amino acid
distribution of an alignment of instances of the motif and entries
from what are known as $PAM$ matrices, which embody patterns of molecular
evolution. In addition, the background
distribution of residues may be modelled more detailed manner, e.g.
using the so-called Dirichlet mixture models.

It is also possible to include position-specific scores for insertion
and deletion of residues, relative to a consensus pattern. When these
are used, the scoring becomes a little more subtle, as the problem is then
quite analogous to pairwise sequence alignment, but with position
dependent scoring parameters for matches, mismatches, insertions and
deletions.

We summarise this section by noting that probability has entered into our
description through the use of frequencies, and scores based on them,
but so far we do not have global statistical models, at least not
ones embodying insertions and deletions, on which we base our estimation
and testing. These are all part of the use of profile HMMs, but first we
introduce HMMs.

\section{Hidden Markov models} \label{section 5}
\setzero\vskip-5mm \hspace{5mm}

Hidden Markov models (HMMs) are processes ${(S_t,O_t), t=1, \ldots, T}$,
where $S_t$ is the hidden state and $O_t$ the observation at time $t$.
Their probabilistic evolution is constrained by the equations
\begin{eqnarray*}
pr(S_t | S_{t-1},O_{t-1},S_{t-2} ,O_{t-2}, \ldots) &=&  pr(S_t | S_{t-1}),\\
pr(O_t | S_{t-1},O_{t-1},S_{t-2} ,O_{t-2}, \ldots) &=&  pr(O_t |
S_t, S_{t-1}).
\end{eqnarray*}
The definitions and basic facts concerning HMMs were laid out in a
series of beautiful papers by L. E. Baum and colleagues around 1970, see
~\cite{durbin:1998} for references.
Much of their formulation has been used almost unchanged  to this day.
Many variants are now used.  For example, the distribution of $O$ may
not depend on previous $S$, or it may also depend on previous $O$ values,
\begin{eqnarray*}
pr(O_t | S_t , S_{t-1} , O_{t-1} ,\ldots )  &=&  pr(O_t | S_t ),\quad
\mbox{or}\\
pr(O_t | S_t , S_{t-1} , O_{t-1} ,\ldots )  &=&  pr(O_t | S_t , S_{t-1}
,O_{t-1}) .
\end{eqnarray*}
Most importantly for us below, the times of $S$ and $O$ may be decoupled,
permitting the observation corresponding to state time $t$ to
be a string whose length and composition depends on $S_t$
(and possibly $S_{t-1}$ and part or all of the
previous observations). This is called a hidden semi-Markov or
generalized hidden Markov model.

Early applications of HMMs were to finance, but these were never published,
to speech recognition, and to modelling ion channels.
In the mid-late 1980s HMMs entered genetics and
molecular biology, where they are now firmly entrenched. One of the
major reasons for the success of HMMs as stochastic models is the fact
that although they are substantial generalizations of Markov chains,
there are elegant dynamic programming algorithms
which permit full likelihood calculations
in many cases of interest. Specifically, there are algorithms which
permit the efficient calculation of a) $pr(sequence|M)$, where $sequence$
is a sequence of observations and $M$ is an HMM; b) the maximum over
$states$ of $pr(states|sequence,M)$, where $states$ is the
unobserved state sequence underlying the observation $sequence$; and c) the
maximum likelihood estimates of parameters in M based on the observation
$sequence$. Step c) is carried out by an iterative procedure which
in the case of independent states was later termed the EM algorithm.

\section{Profile HMMs} \label{section 6}
\setzero\vskip-5mm \hspace{5mm}

\begin{figure}[h]
\centerline{\includegraphics[width=10cm]{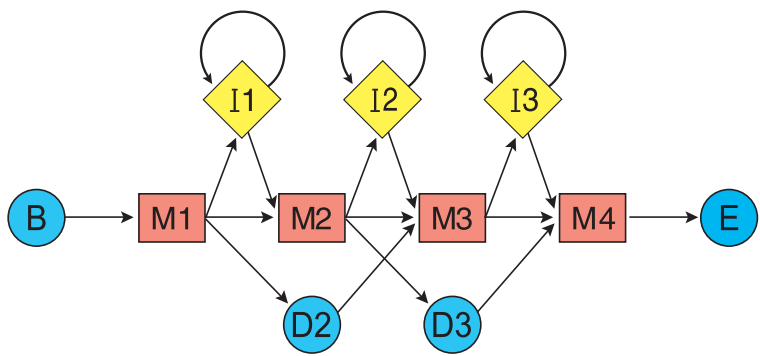}}
  \caption{State space of a simple profile HMM}
  \label{fig3}
\end{figure}

In a landmark paper A. Krogh, D. Haussler and co-workers introduced
profile HMMs into bioinformatics. An illustrative form of their
profile HMM architecture is given in Figure~\ref{fig3}.
There we depict the underlying state space of the hidden Markov chain of
a profile HMM of length 4, with $M$ denoting {\em match} states,
$I$ {\em insert} states and $D$ {\em delete} states, while $B$ and $E$ are
{\em begin} and {\em end} states, respectively.  Encircled states
($D$, $B$ and $E$) do not emit observations, while each of the match and
insert states will have position-specific observation or emission
distributions.  Finally, each arrow will have associated transition
probabilities, with the expectation being that the horizontal transition
probabilities are typically near unity.  This the chain proceeds from
left to right,
and if it remains within match states, its output will be an amino acid
sequence of length 4. Deviation to the insert or delete states will
modify the output accordingly.
The similarity with a direct product of a sequence
of position-specific distributions should be unmistakeable.
The profile HMMs in use now have considerably more features, while sharing
the basic $M, I$ and $D$ architecture.

Why was the introduction of the HMM formalism such an advance? The
answer is simple: it permitted the construction and application of
profiles to be conducted entirely within a formal statistical framework,
and that really helped. Instances of the motif embodied in an HMM could
be identified by calculating $pr(sequence|M)/pr(sequence|B)$ as was done
with profiles, using the algorithm for problem a) in X above. Instances
of the motif could be aligned to the HMM by calculating the most
probable state sequence giving rise to the motif sequence, in essence finding
the most probable sequence of matches, insertions, deletions which align
the given sequence to the others which gave rise to the HMM,  cf. problem b)
above. And finally, the parameters in the HMMs could be estimated
from data comprising known instances of the motif by using maximum
likelihood, an important step for many reasons, one being that it
put insert and delete scores on precisely the same footing as match and
mismatch scores. Although the estimation of HMM parameters
is easiest if the example sequences are properly aligned, the
EM algorithm (problem c) above) does not require aligned sequences. %%%\\

In the years since the introduction of profile HMMs, they have been
become the standard approach to representing motifs and protein
domains. The database Pfam ({\tt http://pfam.wustl.edu})
now has 3,849 hidden Markov models (May 2002) representing
recognized protein or DNA domains or motifs.
Profile HMMs have essentially replaced the use of regular expressions and
the original profiles for searching other databases to find novel
instances of a motif, for finding a motif or domain match to an input
sequence, and for aligning a motif or domain to a an existing
family. There is  considerable evidence that the HMM-based searches are
more powerful than the older profile based ones, though they are slower
computationally, and at times that is an important consideration.

\section{Finding genes in DNA sequence} \label{section 7}
\setzero\vskip-5mm \hspace{5mm}

Identifying genes in DNA sequence is one of the most challenging,
interesting and important problems in bioinformatics today. With so
many genomes being sequenced so rapidly, and the experimental
verification of genes lagging far behind, it is necessary to rely on
computationally derived genes in order to make immediate use of the
sequence.

What is a gene? Most readers will have heard of the famous {\em central
dogma} of molecular biology, in which the hereditary material of an
organism resides in its genome, usually DNA, and where genes are
expressed in a two-stage process: first DNA is {\em transcribed} into
a messenger RNA (mRNA) sequence, and later a processed form of this
sequence is {\em translated}  into an amino acid sequence, i.e. a
protein.
In general the transcribed sequence is longer than the translated portion:
parts called introns (intervening sequence) are removed, leaving exons
(expressed sequence), of which only some are expressed,
 while the rest remain untranslated. The
translated sequence comes in triples called codons, beginning and ending
with a unique start (ATG) and one of three stop (TAA, TAG, TGA) codons.
There are also characteristic intron-exon boundaries called splice donor and
acceptor sites, and a variety of other motifs: promoters, transcription
start sites, polyA sites, branching sites, and so on.\\

All of the foregoing have statistical characterizations, and in
principle they can all help identify genes in long otherwise
unannotated DNA sequence segments. To get an idea of the magnitude
of the task with the human genome, consider the following facts
about human gene sequences ~\cite{nature:2001}: the coding regions
comprise about 1.5\% of the entire genome; the average gene length
is about 27,000 bp (base pair); the average total coding region is
1,340 bp; and the average intron length is about 3,300 bp.
Further, only about 8\% of genes have a single exon. We see that
the information in human genes is very dispersed along the genome,
and that in general the parts of primary interest, the coding
exons, are a relatively
small fraction of the gene, on average about $\frac{1}{20}$.\\

\section{Generalized HMMs for finding genes} \label{section 8}
\setzero\vskip-5mm \hspace{5mm}

\begin{figure}[h]
\centerline{\includegraphics[width=10cm]{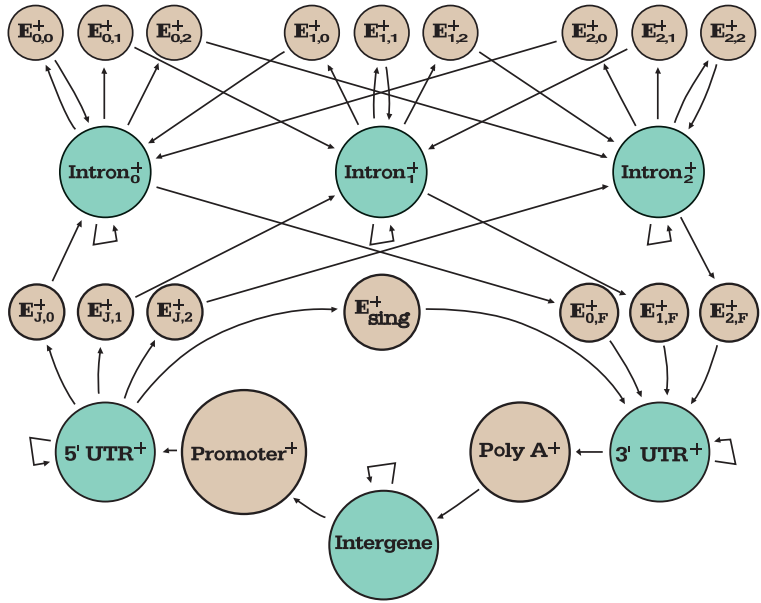}}
  \caption{Forward half of the Genscan GHMM state space}
  \label{fig4}
\end{figure}

The HMMs which are effective in finding genes are the generalized HMMs
(GHMMs) described in section~\ref{section 5} above.
Space does not permit our giving an
adequate description here, so we simply outline the architecture of
Genscan ~\cite{burge:1997}
one of the most widely used human genefinders. States
represent the gene features we mentioned above: exon, intron, and of
course intergenic region, and a variety of other features (promotor,
untranslated region, polyA site, and so on. Output observations embody
state-dependent nucleotide composition, dependence,  and specific signal
features (such as stop codons). In a GHMM the state {\em duration} needs
to be modelled, as well as two other important features of genes in DNA:
the {\em reading frame}, which corresponds to the triples along the mRNA
sequence which are sequentially translated, and the {\em strand}, as DNA is
double stranded, and genes can be on either strand, i.e. they can point
in either direction. These features can be seen in Figure~\ref{fig4},
which was kindly supplied by Lior Pachter.

The output of a GHMM genefinder after processing a genomic segment is
broadly similar to that from a profile HMM after processing an amino
acid sequence: the most probable state sequence given
an observation sequence is a best gene annotation of that sequence, and
a variety of probabilities can be calculated to indicate the support in
the observation sequence for various specific gene features.

%%%
\section{Comparative sequence analysis using HMMs} \label{section 9}
\setzero\vskip-5mm \hspace{5mm}

The large number of sequenced genomes now available, and the observation
that functionally important regions are evolutionarily conserved,
has led to efforts to incorporate conservation into the models and
methods of biological sequence analysis. Pair HMMs were
introduced in \cite{durbin:1998} as a way of including alignment problems
under the HMM framework, and recently \cite{pachter:2002} they were
combined with GHMMs (forming GPHMMs) to carry out alignment and
genefinding with homologous segments of the mouse and human genomes.
Use of the program SLAM on the whole mouse genome
({\tt http://bio.math.berkeley.edu/slam/mouse/}) demonstrated the value of
GPHMMs in this context.

\section{Challenges in biological sequence analysis}
\label{section 10}\setzero\vskip-5mm \hspace{5mm}

The first challenge is to understand the biology well enough to begin
biological sequence analysis. This part will frequently involve
collaborations with biologists. With HMMs, GHMMs and GPHMMs, designing
the underlying architecture, and carrying out the modelling for the
components parts, e.g. for splice sites in genefinding GHMM is perhaps
the next major challenge.  Undoubtedly the hardest and most
important task of all is the implementation: coding up the algorithms
and making it all work with error-prone and incomplete sequence data.
Finally, it is usually a real challenge
to find good data sets for calibrating and evaluating the algorithms,
and for carrying out studies of competing algorithms.

For a recent example of this process, which is a model of its kind, see
~\cite{jarmer:2001}. There an HMM is presented for
the so-called $\sigma^A$ recognition sites, which involve two
DNA motifs separated by a variable number of base pairs. In addition to
the examples mentioned so far, there are many more HMMs in the bioinformatics
literature, see p. 79 of  \cite{durbin:1998} for ones published before 1998.

\section{Closing remarks} \label{section 11}
\setzero\vskip-5mm \hspace{5mm}

In this short survey of biological sequence analysis, I have simply
touched on some of the major ideas. A much more comprehensive
treatment of material covered here can be found in the book
~\cite{durbin:1998}, whose title not coincidentally is the same as
that of this paper. Many important ideas from biological
sequence analysis have not been mentioned here, including molecular
evolution and phylogenetic inference, and the use of stochastic context-free
grammars, a form of generalization of HMMs suited to the analysis of RNA
sequence data.

At this Congress I have talked (and am now writing) on the research of
others, in an area in which my own contributions have been negligible.
I chose to do so upon being honoured by the invitation to speak at this
Congress because I believe this topic -- HMMs -- to be one of the
great success stories of applying mathematics to bioinformatics. In my
view it is the one most worthy of a wider mathematical audience.
I hope that the fact that there are many others better suited than
me to speak on this topic will not prevent readers from appreciating it
and following it up through the bibliography.

I owe what understanding I have of this field to collaborations and
discussions with a number of people, and I would like to acknowledge
them here. Firstly, Tony Wirth, Simon Cawley and Mauro Delorenzi, with
whom I have worked on HHMMs. Next, it has been an honour and pleasure to
observe from close by the development of SLAM, by Simon Cawley,
Lior Pachter and Marina Alexandersson. Finally I'd like to thank
Xiaoyue Zhou and Ken Simpson for their kind help to me when I was
preparing my talk and this paper.

\label{lastpage}

\end{document}